\documentstyle{article}

\newcommand{\ra}{\rightarrow}

\newcommand{\Ra}{\Rightarrow}
\newcommand{\D}{\displaystyle}
\newcommand{\scst}{\scriptscriptstyle}
\newcommand{\scsz}{\scriptsize}

\newcommand{\ot}{\otimes}

\newcommand{\ti}{\times}

\newcommand{\Z}{{\bf Z}}

\newcommand{\ov}{\overline}

\newtheorem{theorem}[section]{Theorem}
{\bfseries}{\itshape}
\newcommand{\proof}{\noindent {\bf Proof.  }}
\newcommand{\qed}{\hfill $\Box$}
\begin{document}
\title{Geometric Objects and Cohomology Operations}
\author{Rocio Gonzalez-Diaz, Pedro Real\\Dept. of Applied Math. I,
University of Seville,  Spain\\
 \{rogodi, real\}@us.es\\
 http://personal.us.es/\{rogodi,real\}
 }
 \date{}

\maketitle

\begin{abstract}
Cohomology operations (including the cohomology ring) of a geometric
 object are finer algebraic invariants than the homology of it.
   In the literature, there exist  various algorithms for computing the homology groups of simplicial complexes
   (\cite{Mun84}, \cite{DE95,ELZ00}, \cite{DG98}), but concerning the algorithmic treatment of
cohomology operations, very little is known. In this paper, we
establish a version of the incremental algorithm for computing
homology given in \cite{ELZ00}, which saves  algebraic
information, allowing us the computation of the cup
product  and  the effective evaluation of the primary and
secondary  cohomology operations on the cohomology of a
finite simplicial complex.
The efficient combinatorial descriptions at cochain level of
 cohomology operations developed
 in \cite{GR99,GR99a} are essential ingredients in our method.
We study the computational complexity
of these processes and a program in Mathematica for cohomology
computations is presented.
\end{abstract}

\section{Introduction}

A simplicial complex is a well--known discrete model of a
geometric object, which consists of a collection of simplices
that fit together in a natural way to form the object. In order
to classify simplicial complexes from a topological point of
view, a first algebraic invariant that can be used is homology,
which in some sense, counts the number of holes of the object.

   We can cite two relevant algorithms for computing homology groups $H_*K$
   of a simplicial complex $K$ in
   ${\bf R}^n$: (1) the classical algorithm
based on reducing certain matrices to their Smith normal form
\cite{Mun84}; (2) the incremental algorithm
\cite{DE95,ELZ00,EZ01}, avoiding the severe computational costs
of the reduction to Smith normal form and consisting of assembling
the complex simplex by simplex and at each step updates the Betti
numbers of the current complex. Starting with the boundary of a
negative simplex, this persistence process finds the cycle which
is destroyed by this simplex through the search, computing in
this way the geometric realization of a homology cycle. It runs
in time at most $O(m^3)$, where $m$ is the number of simplices of
the complex. For simplicial complexes embedding in ${\bf R}^3$,
this complexity is reduced to $O(m)$ in time and space
\cite{DE95}. The algorithm proposed in \cite{DG98} is based on
simulating a thickening of a given complex in ${\bf R}^3$ to a
topological $3$-manifold homotopic to it, and computing the
homology groups of the last one using classical results. The time
and space complexity is linear and this method  also produces
representations of generators of the homology groups.

In general, computing  homology is not enough for determining
whether two geometric objects
    are homeomorphic or not. Finer algebraic invariants  such as
the cohomology (an algebraic dual notion to homology),
    the cup product on cohomology  or cohomology
operations \cite{Spa81}, allow us to topologically distinguish
two geometric objects having isomorphic homology groups. For
example, a torus and the wedge product of a sphere and two
circles have the same homology but the respective cup products on cohomology
are ``essentially'' different.
       Using a field as the coefficient group, for example, ${\bf Z}_2$,
        the cohomology $H^*K$
       of a simplicial complex $K$
       gives us the same topological information as the homology of it. However, the
       additional ring structure on the cohomology determined by the
cup product and cohomology operations cannot directly be produced
from the algorithms previously mentioned for computing the
homology. Roughly speaking, a cohomology operation
$\theta:H^m(-;G) \rightarrow H^n(-;G')$ is a homomorphism that acts on
cohomology ($G$ and $G'$ being groups); relevant examples of cohomology operations are
Steenrod squares, Steenrod reduced powers and Adem secondary
cohomology operations \cite{MT68}. As an example of the strong
constraints that these operations impose on the cohomology of
spaces, we can cite that the use of this machinery is essential
for showing that there do not exist spaces $X$ having cohomology
$H^*(X;{\bf Z})$ a polynomial ring ${\bf Z}[\alpha]$ unless $\alpha$ has
dimension $2$ or $4$.

In this paper, we make use of an
explicit   chain contraction (a special chain equivalence)
connecting the  chain complex $C_* K$, canonically
associated to a simplicial complex $K$ and its homology $H_* K$.
Moreover, from this datum we can derive a cochain contraction
from the cochain complex $C^*K=Hom (C_*K;{\bf Z}_2)$, to the
cohomology $H^* K$.  Using this information, we can compute:
\begin{enumerate}
\item  Geometric realizations of (co)homology generators.
\item The (co)homology class of a (co)cycle in terms of (co)homology generators.
\item The construction of a (co)boundary of a given (co)cycle.
\item  The induced homomorphism at (co)homology level of  a
simplicial map between two complexes.
\item The cup product on cohomology and some primary and secondary cohomology operations.
\end{enumerate}

The first problem is to construct such chain contractions from
$C_* K$ to $H_* K$. In \cite{GR01}, a translation of the
classical matrix algorithm (1) in terms of  chain contractions is
designed. In this paper, we design a version of the incremental
method described in \cite{ELZ00} in terms of chain contractions. The
complexity of our method is  also $O(m^3)$ where $m$ is the number
of simplices of $K$, but our algorithm saves
information which allows us, for example, to compute the following
operations: \begin{enumerate}
 \item The cohomology ring of $K$ in $O(m^5)$.
 \item The Steenrod square operation $Sq^i \alpha_n$ of a cohomology class $\alpha_n$ of
degree $n$ in $O(i^{n-i+1}m)$ (see \cite{GR99a})
 \item The Adem secondary cohomology operation $\Psi_2\alpha_2$ of a cohomology class
 $\alpha_2 \in Ker Sq^2 H^1(K;{\bf \Z_2})$  in $O(m^3)$.
\end{enumerate}
         In fact, the modus operandi for evaluating a mod $2$ cohomology
         operation $\bar{\cal O}: H^m K \rightarrow H^n K$
         on a cohomology class $\alpha_m$ is the
         following:
         \begin{enumerate}
         \item First, given a finite simplicial complex $K$, construct the
         chain contraction from $C^*K$ to $H^*K$
         (denoted $(f^*,g^*,\phi^*):C^*K\Ra H^*K$), using our version of the incremental
         technique.
         \item Evaluate $\bar{\cal O}$ on the cohomology class $\alpha_m$
         using the diagram
          $$\begin{array}{ccc}
          C^{*}K & \stackrel{g^*}\leftarrow & H^*K \\
          {\scst {\cal O}} \downarrow &           & \downarrow \scst{\bar{\cal O}}  \\
           C^*K & \stackrel{f^*}\rightarrow  &H^*K \,,\end{array}$$
           where ${\cal O}: C^*K \rightarrow C^*K$ is a cochain operation associated to
           $\bar{\cal O}$ whose formulation is explicitly given in simplicial
           terms. An efficient combinatorial description
           ${\cal O}$ for $\bar{\cal O}$ being a Steenrod square \cite{GR99,GR99a},
           a Steenrod reduced power \cite{GR99} or
           some Adem secondary cohomology operations \cite{GR01}
           have already been done by the authors. We do not deal
           with this question in this paper, but it is necessary
           to say that the algorithmic approach we give here will only
           be  valid if combinatorial pictures
           of cohomology operations at cochain level are determined.
\end{enumerate}

Let us observe that  in this paper we deal with the general case
of ${\bf R}^n$.
Versions in terms of chain contractions of
 the algorithms given in \cite{DE95} and \cite{DG98}, designed for the special case of ${\bf R}^3$,
would allow us to considerably reduce the computational costs of
the processes.

\section{Homology and Chain Contractions}

In this section, we design a version of the incremental algorithm
of \cite{ELZ00} in terms of  chain contractions. In this way, we construct a  chain contraction
from the chain complex canonically associated to a simplicial
complex $K$, to its homology. Let us observe that passing to
cohomology is not a problem if we use a field as the ground
ring. The resulting cochain contraction from $C^* K$ to $H^* K$
will help us to compute the cup product on cohomology and cohomology operations.

Now, we give a brief summary of concepts and notations. The
terminology
  follows Munkres
 \cite{Mun84}.

Throughout this paper, we consider  ${\bf Z}_2$ is the ground ring and $\mu$ denotes
the product on
${\bf Z}_2$.
A $q$--{\em simplex} $\sigma$ in ${\bf R}^n$(where $q\leq n$) is
 the convex hull of $q+1$ affinely independent points $\{v_0,...,v_q\}$.
 We denote $\sigma=\langle v_0,...,v_q\rangle$.
 The {\em dimension of } $\sigma$ is $|\sigma|=q$.
  A
 $0$--simplex is a vertex, a $1$--simplex is an edge, a
 $2$--simplex is a triangle, a $3$--simplex is a tetrahedron, and
 so on.
       An $i$--{\em face} of $\sigma =\langle v_0,...,v_q\rangle$ ($i<q$) is an
       $i$--simplex whose vertices are in the set
  $\{ v_0,...,v_q\}$. The $(q-1)$--faces of $\sigma$ are  called the {\em facets} of $\sigma$.   A
simplex is {\em shared} if it is a face of more than one simplex.
Otherwise, the simplex is {\em free} if it belongs to one
higher--dimensional simplex, and {\em maximal} if it does not
belong to any.
          A {\em simplicial complex} $K$ is a collection of simplices such that:
\begin{itemize}
\item If $\tau$ is a face of $\sigma\in K$, then $\tau\in K$.
\item If $\sigma', \sigma\in K$, then $\sigma'\cap\sigma\in K$ or $\sigma'\cap\sigma=\emptyset$.
\end{itemize}
Let us notice that $K$ can be given by the set of its maximal simplices. The {\em dimension of }
$K$ is $dim K=max \{|\sigma|:\; \sigma\in K\}$.
In this paper, all the simplices have finite dimension and all
the simplicial complexes are finite collections. The set of all the
$q$--simplices of $K$ is denoted by $K^{(q)}$.
If $L$ is a
subcollection of $K$ that contains all faces of its elements,
then $L$ is a simplicial complex in its own right; it is called a
{\em subcomplex} of $K$.
Let $K$ and $K'$ be two
simplicial complexes. A  map $f:K^{(0)}\rightarrow K'^{(0)}$ such that
whenever $\langle v_0,...,v_q\rangle \in K$ then $f(v_0),...,f(v_n)$ are
vertices of a simplex of $K'$, is called a {\em vertex map}.

Algebraic Topology is the study of algebraic objects attached to
topological spaces; the algebraic invariants reflect some of the
topological structure of the spaces.

The {\em chain complex} $C_{*} K$ associated to a simplicial
complex $K$ is a family $\{C_q K,\partial_q\}_{q\geq 0}$ defined
in each dimension  $q$ by:
\begin{itemize}
\item $C_q K$ is the free abelian group generated by the  $q$--simplices of $K$.
An  element $a=\sigma_1+\cdots+\sigma_m$ of
$C_q K$ ($\sigma_i\in K^{(q)}$) is called a
$q$--{\em chain}.
\item $\partial_q:\,C_q K\rightarrow C_{q-1} K$ called the {\em boundary operator} is given by
$$\displaystyle
\partial_q \langle v_0,...,v_q\rangle=\sum_{i=0}^{q}\langle v_0,...,\hat{v}_i,...,v_q\rangle$$
where $\langle v_0,...,v_q\rangle$ is a $q$--simplex
of $K$ and
 the hat means that $v_i$ is omitted.
 By
 linearity,  $\partial_q$ can be extended to
 $C_q K$, where it is a homomorphism.
 \end{itemize}
 A $q$--chain $a$ is called a $q$--{\em cycle} if $\partial a=0$.
  If $a= \partial b$ for some $b\in C_{q+1} K$ then $a$
 is called a $q$--{\em boundary}.
  We denote the groups of $q$--cycles
 and $q$--boundaries by $Z_q K$ and $B_q K$ respectively, and
 define $Z_0 K=C_0 K$.
   Since $B_q K\subseteq Z_q K$, we can
 define the  {\em  $q$th homology group} to be the quotient group
 $Z_q K/B_q K$, denoted by $H_q K$.
  Given that elements of this
 group are cosets of the form $a+B_q K$, where $a\in
 Z_q K$, we say that the coset $a+B_q K$, denoted by $[a]$, is the {\em homology
 class} in $H_q K$ determined by $a$ or $a$ is a
  {\em representative cycle} of $[a]$.
  Let $K$ and $L$ be two simplicial complexes.
   A {\em chain map} $f:C_* K\rightarrow C_* L$ is a family of homomorphisms
  $$\{f_q:C_q K\rightarrow C_q L\}_{q\geq 0}$$
  such that $\partial_q f_q=f_{q-1}\partial_q$ for all $q$.
Observe that for every vertex map
 $f:K^{(0)}\rightarrow L^{(0)}$, we can obtain the corresponding chain map
 $f_{\#}:C_* K\rightarrow C_* L$ such that
$$\displaystyle f_{\#}\langle v_0,...,v_q\rangle =\left\{
\begin{array}{cl}
\langle f(v_0),...,f(v_q)\rangle\quad & \mbox{ if $f(v_i)$ distinct}\\
0&\mbox{ otherwise}
\end{array}
\right.$$

Let  $h$ and $k$ be two chain maps from $C_* K$ to $C_* L$. A
{\em chain homotopy} from $h$ to $k$ is a family of homomorphisms
$$\{\phi_q: C_q K\rightarrow C_{q+1} L\}_{q\geq 0}$$ such that
$\partial_{q+1} \phi_q+\phi_{q-1}\partial_q=h_q+k_q$. We write
$h\sim k$ if a chain homotopy between $h$ and $k$ exists. Two
chain complexes $C_* K$ and $C_* L$ are {\em chain equivalent} if
there exist two chain maps $f:C_* K\rightarrow C_* L$ and $g:C_*
L\rightarrow C_* K$ such that $$ fg\sim 1_{\scst C_* L}\qquad
\mbox{ and } \qquad gf\sim 1_{\scst C_* K}\,.$$ Observe that, in
this case, $\phi_q: C_q K\rightarrow C_{q+1} K$ for  all $q\geq
0$. A {\em chain contraction} \cite{EM52} from $C_* K$ to $C_* L$
is a chain equivalence
 such that
$$fg=1_{\scst C_* L}\qquad\mbox{ and } \qquad gf\sim 1_{\scst C_* K}
\quad \mbox{(that is, $1_{\scst C_* K}+gf=\partial \phi+\phi\partial$)}$$
and $\phi$ has the following ``annihilation"
properties: $f\phi=0,\; \phi g=0$  and $\phi\phi=0\,.$
We denote such chain contraction as $(f,g,\phi): C_*K\Ra C_*L$.
Observe that if a   chain contraction from $C_* K$
to $C_* L$ exists then  $L$ has fewer or the same number of
simplices than $K$.
Now, we show some examples of contractions.
\begin{itemize}
\item[(a)] Edge Contractions.\label{edge}

Conditions under which edge contractions are homeomorphisms appear in \cite{DEGN99}.
Here, we show
one condition under  which edge contractions become, at algebraic level, chain contractions.

Let $K$ be a simplicial complex and
 $\tau=\langle a,b\rangle$ an edge in
$K$. An {\em edge contraction} is  given by the vertex map $\;f:
K^{(0)}\rightarrow L^{(0)}=K^{(0)}-\{a,b\}\cup \{c\}\;$ where $f(a)=f(b)=c$,
and $f(v)=v$ for all $v\neq a,b$.

Let $B$ be a subset of $K$ that is not necessarily a subcomplex.
Define
$$\ov{B}=\{\sigma'\in K:\;\; \sigma'\leq \sigma\in B\},\qquad
St\,B=\{\sigma\in K:\;\;\sigma\geq
\sigma'\in B\}\,,$$
$$ Lk\,B=\overline{St\,B}-St\,\ov{B}\,,$$
where $\sigma'< \sigma$ means that $\sigma'$ is a face of $\sigma$.

If $Lk\,a\,\cap\, Lk\,b= Lk\,\tau,$
  then a
 chain contraction $(f_{\#},g,\phi)$ from $C_*K$ to $C_*L$ is defined as follows:
\begin{itemize}
\item $f_{\#}$ is the chain map induced by the vertex map $f$.
\item $g:C_*L\rightarrow C_*K$ is such that
$$\begin{array}{l}
g \tau=\tau \qquad \forall\tau\not\in St\,c,\\
g\langle c\rangle=\langle a\rangle,\\
g (\omega\cup \langle c\rangle)=\left\{
\begin{array}{cll}
   \omega\cup \langle a\rangle  &\mbox{ if }& \omega \in Lk\,a,\\
  \omega\cup  \langle b\rangle+ \bar\omega\cup \langle a,b\rangle
&\mbox{ if } &\omega \in Lk\,b - Lk\,\tau \\
&& \bar{\omega}\in Lk\,\tau\mbox{ and } \bar{\omega}<\omega,\\
 \omega\cup \langle b\rangle
&\mbox{ if }&\omega \in Lk\,b - Lk\,\tau,\\
&&\not\exists \,\bar{\omega}<\omega\mbox{ and }\bar{\omega}\in Lk\,\tau.\end{array}\right.
 \end{array}$$
\item  $\phi: C_*K\rightarrow C_{*+1}K$ is  given by
$$\begin{array}{ll}
\phi \langle v_0,...,v_q,b\rangle=\langle v_0,...,v_q,a,b\rangle &\mbox{ if
$\langle v_0,...,v_q\rangle \in Lk\,\tau$}
\end{array}$$
and $\phi\tau=0$ otherwise.
\end{itemize}

\item[(b)] Simplicial Collapses.\label{collapse}

 Suppose $K$ is a simplicial complex, $\sigma\in K$ is a maximal
$q$--simplex
  and $\sigma'$ is a free $(q-1)$--face of $\sigma$.
 Then,  $K$ {\em simplicially collapses} onto $K-\{\sigma, \sigma'\}$.
 More generally, a {\em simplicial collapse} is any sequence of
 such operations.
  A {\em thinned} simplicial complex
  $M_{\mbox{\scsz scol}}(K)$ is a subcomplex of $K$ with the condition
that
  all the faces
 of the maximal simplices of $M_{\mbox{\scsz scol}}(K)$ are shared.
 Then, it is
 obvious that it is no longer possible to collapse.
 There is an explicit chain contraction from $C_*K$ onto
 $C_*(M_{\mbox{\scsz scol}}K)$ \cite{For99}.
 The following algorithm computes $M_{\mbox{\scsz scol}}K$ and
 the chain contraction from $C_*K$ onto
 $C_*(M_{\mbox{\scsz scol}}K)$.
 Suppose that $K$ is given by the set of its maximal simplices.
 \begin{center}
 \begin{tabbing}
 {\tt Initially, $M_{\mbox{\scsz scol}}K=K$, $\;\;\phi\tau=0$,
 $\;f\tau=g\tau=\tau$ for each $\tau\in K$}.\\
 {\tt While} \= {\tt there exists a maximal simplex $\sigma$ with a
free face $\sigma'$ do}\\
  \> {\tt $M_{\mbox{\scsz scol}}K=M_{\mbox{\scsz scol}}K-\{\sigma,
\sigma'\},$}\\
\> {\tt $\phi\sigma'=\sigma$, $f\sigma'=\sigma'+\partial\sigma$ and $\;f\sigma=0$}\\
 {\tt End while}
 \end{tabbing}
 \end{center}

\item[(c)] Contraction to a Vertex.\label{vertex}

Let $\sigma=\langle v_0,\dots,v_q\rangle$ be a simplex and let $K[\sigma]$ be the simplicial complex whose maximal simplex is
$\sigma$. It is obvious that we can obtain a chain contraction from $C_*K[\sigma]$ to $\langle v_0\rangle$
using simplicial collapses. But now, we show another contraction from $C_*K[\sigma]$ to
$\langle v_0\rangle$ determining the acyclicity of the simplex $\sigma$. This last chain contraction
is the key for constructing another one from any simplicial complex to its homology as we will see
in the following section.
We define $(f_{\sigma},g_{\sigma},\phi_{\sigma}): C_*K[\sigma]\Ra \langle v_0\rangle$ as follows:
\begin{eqnarray*}
&&f_{\sigma}\langle v_i\rangle=\langle v_0\rangle\qquad 0\leq i\leq q\,,\quad\mbox{ and }\quad
f_{\sigma}(\tau)=0\quad \mbox{otherwise};\\
&&g_{\sigma}\langle v_0\rangle=\langle v_0\rangle;\\
&&\phi_{\sigma}\langle v_0,v_{j_1},\dots,v_{j_n}\rangle=0\quad\mbox{and}\quad
\phi_{\sigma}\langle v_{j_1},\dots,v_{j_n}\rangle=\langle v_0,v_{j_1},\dots,v_{j_n}\rangle
\end{eqnarray*}
where $1\leq j_1<\cdots<j_n\leq q$.

Let us observe that in this case $\langle v_0\rangle$ represents the unique class of homology in
$H_*K[\sigma]$.

\end{itemize}

\subsection{Incremental Homology Algorithm and  Chain Contractions}

Our algorithm  for computing a chain contraction from the chain
complex of a simplicial complex $K$ to its homology is based on
the incremental algorithm  for computing the persistence of the
Betti numbers developed in \cite{ELZ00}.

The input of our algorithm implemented in Mathematica is
the sorted set of all the simplices, $K=\{\sigma_1,\dots,\sigma_m\}$,
with the property that any subset of it, $\{\sigma_1,\dots,\sigma_i\}$, $i\leq m$, is a simplicial
complex itself.
The output $\ell=$ {\tt contraction[$K$]} is a list of sorted lists.
Each sorted list has three elements.
The first one
is a simplex $\sigma$ of $K$, the second one is the image of $\sigma$ under $f$ and the third
one consists of
the image  of $\sigma$ under $\phi$. We omit in the list the simplices such that the image of them
 are null
under $f$ and $\phi$.
 In general, a class of  homology $\alpha$ is represented by a simplex $\tau$,
 so in order to obtain the image of $\alpha$ under $g$, we only have to compute $a=
\tau+\phi \partial \tau$. Moreover, $a$ will be a representative
cycle of $\alpha$.

Now, let us suppose we have
constructed the list $\ell=${\tt contraction[$L$]}
for $L=\{\sigma_1,\dots,\sigma_{i-1}\}$, $i\leq m$ (if $L=\emptyset$,
we assume $\ell=\emptyset$). We construct
{\tt contraction[$\{\sigma_1,\dots,\sigma_{i}\}$]}
as follows:
\begin{center}
\begin{tabbing}
$\mbox{ }\mbox{ }$\= {\tt If $\mbox{ }\mbox{ }$}\= {\tt $f$[$\partial\sigma_i,\ell$]$=0$ then,}\\
\>\> $\ell\cup \{(\sigma_i,\sigma_i,\phi \sigma_i)\}$,\\
\>{\tt Else}\\
\>\> {\tt Replace} \={\tt [} \=  {\tt Replace} \= {\tt [} \= $\ell$,\\
\>\>\>\>\>\> {\tt Solve[$f$[$\partial\sigma_i,\ell$]=$0$]}\\
\>\>\>\>\>{\tt ],}\\
\>\>\>\>{\tt Solve[$\phi$[$\partial\sigma_i,\ell$]=$\sigma_i$]}\\
\>\>\>{\tt ]}\\
\>{\tt End if}
\end{tabbing}
\end{center}
where,  for a simplex $\tau$, $f[\tau,\ell]$ and $\phi[\tau,\ell]$ are, respectively,
the second and the third element of the list of $\ell$ that has $\tau$ as the first element.
If this list does not exist, then {\tt $f$[$\tau,\ell$]$=0$} and {\tt $\phi$[$\tau,\ell$]$=0$}.
Now, let us explain what {\tt contraction[$\{\sigma_1,\dots,\sigma_{i}\}$]} computes.
If {\tt $f$[$\partial \sigma_i,\ell$]$=0$} then $\sigma_i$ ``creates a
cycle", so in fact,
$\sigma_i$ is a new generator of homology.
Otherwise, {\tt $f$[$\partial \sigma_i,\ell$]} is a sum of elements of the form
$\sum_{\sigma_j\in N\subset L} \sigma_j$. The idea of this last case is that
$\sigma_i$ destroys the cycle generated by $\partial \sigma_i$ in $L$. Therefore,
we impose {\tt $f$[$\partial \sigma_i,\ell$]$=0$} and
{\tt $\phi$[$\partial \sigma_i,\ell$]$=\sigma_i$}.
We replace these relations in $\ell$ with the commands Replace and Solve.

At the end of the algorithm, all the elements of the form $\phi\tau$ are replaced by zero.
For obtaining the morphism $g$ and the representative cycles of the homology classes
of $K$,
we compute $\tau+\phi \partial\tau$ for each simplex $\tau$ (the generators of homology)
 satisfying that {\tt $f$[$\tau,\ell$]$=\tau$} in the list
 $\ell=${\tt contraction}[$K$]. We create a new list of sorted lists, called
{\tt representativeCycles}[$K$] such that in each sorted list
the first element is a generator of homology, $\tau$, and the second element is its image under
$g$, $\tau+\phi\partial\tau$.
Observe that this last chain is, in fact,  a cycle:
\begin{eqnarray*}
\partial (\tau+\phi \partial\tau)&=&\partial \tau+\partial \phi \partial\tau\\
&=&\partial \tau+(gf-1-\phi\partial)\partial\tau \\
&=&gf \partial \tau-\phi\partial\partial\tau \qquad \mbox{[ since $\partial\partial\tau=0$,
then ]}\\
&=&gf \partial \tau \qquad \mbox{[ since, by construction,  $f\partial\tau=0$, then ]}\\
&=&0.
\end{eqnarray*}
It is easy to check that $(f,g,\phi)$ is, in fact, a chain contraction from $C_*K$ to $H_*K$.
Observe that given a cycle $a$, if $fa=0$
then $a$ is  also a boundary.
In order to compute a chain $a'$ such
 that $a=\partial a'$, we can use the relation
 $$a-g f a=\phi\partial a+\partial\phi a\,.$$
Since  $\partial a=0$ and $fa=0$, we have $a=\partial\phi a\,.$

\begin{theorem}
{\em The complexity of our algorithm for computing the homology of a finite simplicial
complex $K$
and a chain contraction from  $C_*K$ on $H_*K$  is
$O(m^3)$, where $m$ is the number of simplices of $K$.}
\end{theorem}

\proof

Let $K=\{\sigma_1,\dots,\sigma_m\}$ and $d=dim\, K$. Suppose that we have computed
$\ell=${\tt contraction[$\{\sigma_1,\dots,$ $\sigma_{i-1}\}$]}.
In the worst case,  we have to solve {\tt $f$[$\partial \sigma_i,\ell$]} $=0$ and
{\tt $\phi$[$\partial \sigma_i,\ell$]$=\sigma_i$}. Observe that
the number of
simplices involved in $\partial \sigma_i$ is less or equal than the dimension of  $\sigma_i$
which is at most $d$
and then, the number of
simplices involved in the formulas of {\tt $f$[$\partial \sigma_i,\ell$]} and
{\tt $\phi$[$\partial \sigma_i,\ell$]} is
$O(d m)=O(m)$. Since we have to solve the equations and replace the solution in $\ell$,
the total cost of these operations is
$O(m^2)$.
Moreover,
for obtaining the representative cycles, we have to compute $\tau+\phi\partial\tau$ for every
generator of homology. The cost of this is also $O(m^2)$.
 Therefore, the total algorithm runs in time at most $O(m^3)$.
\qed

\section{Cohomology and Cohomology Operations}

One reason in order to use the cohomology for distinguishing spaces instead of homology, is that the cohomology
has additional structures, such as  cup product and cohomology operations. If two spaces
have isomorphic (co)homology groups but the
behaviour of the ring structure or cohomology operations is different, then
they are not homeomorphic. In this section we explain how we can compute the cup product and
cohomology operations starting
from a chain contraction from  an algebraic object to its homology. We first need to define more
concepts.

  The {\em cochain complex} associated to $K$, denoted by
 $C^* K$, is the family
 $$\{C^q K,\delta^{q}\}_{q\geq 0}\,,$$
 defined in each dimension $q$ by:
 \begin{itemize}
\item  The group $C^q K=Hom (C_q K;{\bf Z}_2)$=$\{c:\,C_q K\rightarrow {\bf Z}_2\;$ is
 a homomorphism$\}$.
 \item  The homomorphism $\delta^q:C^q K\ra
 C^{q+1}K$ called the {\em coboundary operator}  given by
 $$\delta^q
  c\, a=c \,\partial_{q+1}a $$ where $c\in C^q K$ and
 $a\in C_{q+1} K$.
 \end{itemize}
 The elements of $C^q K$ are called $q$--{\em cochains}.
 Observe that a $q$--cochain can be defined  on $K^{(q)}$
 and it is naturally extended by linearity on $C_q K$.
    $Z^q K$ and $B^q K$ are the kernel of $\delta^q$ and the
 image of $\delta^{q-1}$, respectively.
  The elements in $Z^q K$
 are called $q$--{\em cocycles} and those in $B^q K$ are called
 $q$--{\em coboundaries}.
  The {\em $q$th cohomology group}
 $$H^q K=Z^q K/B^q K$$
 can be defined for each integer
 $q$.
 Take into account that since the ground ring is a field, the homology and cohomology of $K$
are isomorphic. Moreover, given a generator of homology, $\alpha$, of dimension $q$,
we can define the corresponding generator of cohomology $\alpha^*: H_qK\ra {\bf Z}_2$ such as
$$\alpha^* \alpha=1\qquad \mbox{ and }\qquad \alpha^*\beta=0\quad\mbox{ for }
\alpha\neq \beta\in H_qK\,.$$
One can also define the dual concept of chain maps and
chain contractions, in the obvious way.
Furthermore, starting from a  chain contraction $(f,g,\phi)$ from $C_* K$ to $H_* K$, we
construct a cochain contraction $(f^*,g^*,\phi^*)$ from $C^* K$ to $H^* K$ as
follows. Let $c\in C^* K$ and $\alpha^*\in H^* K$. Define $f^* c=c\,
g$, $g^* \alpha^*=\alpha^* f$ and $\phi^* c=c\,\phi$.

The cohomology of  $K$ is a ring with the
{\em cup product} $$\smile: \; H^iK\ot H^jK\ra H^{i+j}K$$
defined at a cocycle level by $(c\smile c')\sigma=\mu
(c\langle v_{0},\dots,v_{i}\rangle \ot c'\langle v_i,\dots,v_{i+j}\rangle)$, where $c$ and
$c'$ are an $i$--cocycle and a $j$--cocycle, respectively, and
$\sigma=\langle v_0,\dots,v_{i+j}\rangle\in K^{(i+j)}$ is such that
$v_0<\cdots < v_{i+j}$.
Using the chain contraction $(f,g,\phi)$ from $C_*K$ to $H_*K$,
we can compute the cohomology ring of $K$ in
the following way:
\begin{center}
\begin{tabbing}
{\tt Take $\alpha^*$ and $\beta^*$, cohomology classes of $K$}\\
{\tt For }\={\tt every $\gamma\in H_{i+j}K$}\\
\>{\tt compute $((\alpha^* f)\smile (\beta^* f))g\gamma$}\\
{\tt End for} \end{tabbing}
\end{center}
Notice that the resulting cohomology class  is
 determined by the cocycle $c=(\alpha^* f)\smile (\beta^* f)$.

In order to compute a cohomology operation $\bar{\cal O}: H^*K\ra
H^{*+i}K$, on one hand, we need to compute
 {\tt contraction[$K$]} in order to obtain a chain contraction $(f,g,\phi)$ from $C_*K$ to
its homology and, on the other hand, we need a simplicial
version ${\cal O}:C^*K\ra C^{*+i}K$ of $\bar{\cal O}$. Therefore, for obtaining
$\bar{\cal O}(\alpha^*)$,
where $\alpha^*\in H^*K$, we only need to compute ${\cal O}(\alpha^*f)g$ (for more details, see
\cite{GR01}).
For example, from the combinatorial formulae of Steenrod squares given in \cite{Ste47,SE62},
 $$Sq^i:H^*K\rightarrow  H^{*+i}K\,,$$
 for calculating the cohomology class $Sq^i(\alpha^*)$ with $\alpha^*$ in $H^qK$, we
only have to compute $Sq^i(\alpha^* f)g$.
More concretely, at cochain
level, $Sq^ic=c\smile_{q-i}c$ mod $2$. Moreover, given a $p$--cochain $c$ and a
$q$--cochain $c'$,
$c\smile_n c'$ is a $(p+q-n)$--cochain defined by
$$\D (c\smile_n c')\sigma=  \sum_{0\leq i_0<\cdots<i_n\leq p+q-n}\mu(c
(\cup_{\mbox{\scsz $j$ even}} z^j)
\ot c'(\cup_{\mbox{\scsz $j$ odd}} z^j))$$
where $\sigma=\langle v_0,\dots,v_{p+q-i}\rangle$, $v_0<\cdots<v_{p+q-i}$;
$z^0=\langle v_0,\dots,v_{i_0}\rangle$,
$z^j=\langle v_{i_{j-1}},\dots, v_{i_j}\rangle$, for $1\leq j\leq n$,
and $z^{n+1}=\langle v_{i_n},\dots, v_{p+q-n}\rangle$.
Finally,
we can express Steerond squares in
a matrix form due to the fact that these cohomology operations
are homomorphisms. The process of diagonalization of
such matrices can give us  detailed information about the kernel
and image of these cohomology operations. This information will be very useful in the next section in
order to compute Adem secondary cohomology operations.

\section{Adem Secondary Cohomology Operations}

For attacking the computation of secondary cohomology operations,
we will see in this section that the homotopy operator $\phi$
of the chain contraction $(f,g,\phi)$ from $C_*K$ to  the homology of $K$, is essential.

First of all,
we will need the following mod $2$ relation \cite{Ste47}:
\begin{eqnarray}
\label{cup}
\delta(c\smile_n c')=c\smile_{n-1}c'+c'\smile_{n-1}c+\delta c\smile_n c'+c\smile_n \delta c'
\end{eqnarray}
where $c$ and $c'$ are two cochains.
Now, we shall indicate how Adem secondary cohomology
operations
\begin{eqnarray*}
\Psi_q: N^qK\ra H^{q+3}(K;{\bf Z}_2)/Sq^2 H^{q+1}(K;{\bf Z}),\qquad q\geq 2
\end{eqnarray*}
can be constructed (see \cite{Ade52,Ade58}). $N^qK$ denotes the
kernel of $Sq^2: H^q(K;{\bf Z})\ra H^{q+2}(K;{\bf Z}_2)$.
These operations appear using the known relation:
 $$Sq^2Sq^2\alpha+Sq^3Sq^1\alpha=0$$
for any $\alpha\in H^*(K;{\bf Z})$. For this particular relation
there exist cochain mappings
$$E_{j}:C^*(K\ti K\ti K\ti K)\ra C^{*-j}K$$ such that mod $2$
$$(c\smile_{q-2}c)\smile_q (c\smile_{q-2}c)
+(c\smile_{q-1}c)\smile_{q-2} (c\smile_{q-1}c)=\delta
E_{3q-3} c^4\,,$$ where  $c$ is a $q$--cochain with integer coefficients.
Making use of the relation (\ref{cup}) we have that mod $2$
$$(c\smile_{q-2}c)\smile_q (c\smile_{q-2}c)=\delta (b\smile_q\delta b+b\smile_{q-1}b)$$
$$(c\smile_{q-1}c)\smile_{q-2} (c\smile_{q-1}c)=
\delta(\eta\smile_{q-2}\delta\eta+\eta\smile_{q-3}\eta)$$
where $b$ is a $(q+1)$--cochain such that
$c\smile_{q-2}c=\delta b$ and $\eta=\frac{1}{2}(c\smile_q
c+c)$.
Therefore
$$w=\left\{\begin{array}{l}
E_{3q-3}c^4
+b\smile_{q-1} b+b\smile_{q}\delta b+\eta\smile_{q-2}\delta\eta+\eta\smile_{q-3}\eta,\qquad q>2\\\\
E_{3}c^4
+b\smile_1 b+b\smile_{2}\delta b+\eta\smile\delta\eta,\qquad q=2
\end{array}\right.$$
is a mod $2$ cocycle. If
 $c$ is a representative $q$--cocycle of a cohomology class $\alpha\in N^q K$ with integer
 coefficients then,
$$\Psi_q \alpha=[w]+Sq^2 H^{q+1} K\,.$$

Now, suppose ${\bf Z}_2$ is the ground ring and suppose we have computed
the contraction $(f,g,\phi)$
from $C_*K$ to $H_*K$,
$\ell=${\tt contraction[$K$]}.
Then, the cochain $b$ is $\phi^*(c\smile_{q-2}c)=(c\smile_{q-2}c)\phi$.
Observe that for computing $\Psi_q \alpha^*$, $\alpha^*\in H^*K$,
we need to have a combinatorial expression of the
morphism $E_{3q-3}$.
A method for
obtaining ``economical" combinatorial formulae for $E_{3q-3}$ is given in \cite{Gon00}. For example,
 \begin{eqnarray*}
(E_3 c^4 )\sigma
&=& \mu( c\langle v_0,v_2,v_3\rangle \ot c(v_0,v_1,v_2\rangle
\ot c\langle v_3,v_4,v_5\rangle\ot c\langle v_2,v_3,v_5\rangle\\
    &&+c\langle v_0,v_4,v_5\rangle\ot c\langle v_3,v_4,v_5\rangle
    \ot c\langle v_0,v_1,v_2\rangle\ot c\langle v_0,v_1,v_2\rangle\\
    &&+c\langle v_0,v_1,v_5\rangle\ot c\langle v_3,v_4,v_5\rangle
    \ot c\langle v_1,v_2,v_3\rangle\ot c\langle v_1,v_2,v_3\rangle\\
    &&+c\langle v_0,v_1,v_2\rangle\ot c\langle v_2,v_4,v_5\rangle
    \ot c\langle v_2,v_3,v_4\rangle\ot c\langle v_2,v_3,v_4\rangle\\
    &&+c\langle v_0,v_1,v_2\rangle\ot c\langle v_2,v_3,v_5\rangle
    \ot c\langle v_3,v_4,v_5\rangle\ot c\langle v_3,v_4,v_5\rangle)\,,
\end{eqnarray*}
where $c$ is a $2$--cochain and   $\sigma=\langle v_0,v_1,v_2,v_3,v_4,v_5\rangle$ is a
$5$--simplex such that $v_0<v_1<v_2<v_3<v_4<v_5$. Therefore, the steps for computing
$\Psi_q$ are the following:
\begin{itemize}
\item[1.] Take $\alpha^*\in N^qK$ making use of the diagonalization of the matrix of
$Sq^2H^qK$.
\item[2.] Compute $c=\alpha^*f$.
\item [3.] Compute $\;b=(c\smile_{q-2}c)\phi$,
$\;\eta=\frac{1}{2}(c\smile_q c+c)$,
$\; b\smile_{q-1} b$, $\; b\smile_q\delta b$, $\;\eta\smile_{q-3}\eta$,
$\; \eta\smile_{q-2}\delta\eta\;$ and  $\;E_{3q-3}c^4$.
\item[2.] Compute $wg$.
\end{itemize}
Let us explain with more detail the first step.  In our implementation in Mathematica,
the command {\tt hclass[$\ell,q$]} computes
the list of all the cohomology classes of $K$
in dimension $q$.
We compute $Sq^2\alpha^*$ for each $\alpha^*\in$ {\tt hclass[$\ell,q$]} and we write
the result as a
vector {\tt sq2[ $\ell,\alpha^*$]}
of $0's$ and $1's$ such that $$\mbox{$Sq^2\alpha^*=${\tt sq2[ $\ell,\alpha^*$]}. {\tt hclass[$\ell,q+2$]} }.$$
Then, we construct the matrix corresponding to $Sq^2H^qK$ with the
command
$$\begin{array}{l}
\mbox{\tt matrixSq2[$\ell,q$]}\\
\qquad\qquad\mbox{ $=${\tt Table[sq2[}$\ell$,
{\tt hclass[$\ell,q$][[$i$]]]},$\{i,1,${\tt Length[hclass[$\ell,q$]]}\}}\end{array}$$
After this, we compute
$$\mbox{{\tt NullSpace[matrixSq2[$\ell,q$],
 Modulus$\rightarrow 2$]}. {\tt hclass[$\ell,q$]}}$$ in order to obtain
a base of $N^qK$.

An example of the computation of Adem secondary cohomology operation using
our algorithm is the following.
Let $K$ be a simplicial complex whose set of maximal simplices is
\begin{eqnarray*}
\{&\langle 1, 3, 7\rangle, \langle 3, 4, 7\rangle, \langle 1, 4, 7\rangle,
\langle 1, 2, 8\rangle,
\langle 2, 3, 8\rangle, \langle 1, 3, 8\rangle,&\\
    &\langle 4, 5, 9\rangle, \langle 4, 6, 9\rangle, 
\langle 5, 6, 9\rangle,
     \langle 3, 4, 10\rangle,
\langle 3, 6, 10\rangle, \langle 4, 6, 10\rangle,&\\
     &\langle 1, 2, 3, 4, 5, 6\rangle, \langle 1, 2, 3, 4, 5, 11\rangle,
     \langle 1, 2, 3, 4, 6, 11\rangle,&\\&
  \langle 1, 2, 3, 5, 6, 11\rangle, \langle 1, 2, 4, 5, 6, 11\rangle,
  \langle 1, 3, 4, 5, 6, 11\rangle,
    \langle 2, 3, 4, 5, 6, 11\rangle&\}
   \end{eqnarray*}
     We first compute the chain contraction to the homology:
\begin{eqnarray*}
\{&\{\langle 1\rangle, \langle 1\rangle, 0\},
\{\langle 2\rangle, \langle 1\rangle,\langle 1, 2\rangle\},
\{\langle 3\rangle,    \langle 1\rangle, \langle 1, 3\rangle\},&\\
&\{\langle 4\rangle, \langle 1\rangle, \langle 1, 3\rangle + \langle 3, 4\rangle\},
\{\langle 5\rangle, \langle 1\rangle, \langle 1, 3\rangle +\langle 3, 4\rangle + \langle 4, 5\rangle\},&\\
&\{\langle 6\rangle, \langle 1\rangle,  \langle 1, 3\rangle + \langle 3, 4\rangle + \langle 4, 6\rangle\},
\{\langle 7\rangle, \langle 1\rangle, \langle 1, 7\rangle\},
\{\langle 8\rangle, \langle 1\rangle,  \langle 1, 8\rangle\},&\\
&\{\langle 9\rangle,    \langle 1\rangle,\langle 1, 3\rangle + \langle 3, 4\rangle + \langle 4, 9\rangle\},
\{\langle 10\rangle, \langle 1\rangle,    \langle 1, 3\rangle + \langle 3, 10\rangle\},&\\
&\{\langle 11\rangle, \langle 1\rangle,   \langle 1, 11\rangle\},
\{\langle 1, 4\rangle, 0,   \langle 1, 3, 7\rangle + \langle 1, 4, 7\rangle+ \langle 3, 4, 7\rangle\},&\\
&\{\langle 1, 5\rangle, 0,    \langle 1, 3, 7\rangle + \langle 1, 4, 5\rangle
+ \langle 1, 4, 7\rangle
+\langle 3, 4, 7\rangle\},&\\
&\{\langle 1, 6\rangle, 0,    \langle 1, 3, 6\rangle + \langle 3, 4, 10\rangle + \langle 3, 6, 10\rangle
+ \langle 4, 6, 10\rangle\},&\\
&\{\langle 2, 3\rangle, 0,    \langle 1, 2, 8\rangle + \langle 1, 3, 8\rangle
+ \langle 2, 3, 8\rangle\},&\\
&\{\langle 2, 4\rangle, 0,    \langle 1, 2, 4\rangle + \langle 1, 3, 7\rangle
+ \langle 1, 4, 7\rangle +
\langle 3, 4, 7\rangle\},&\\
&\{\langle 2, 5\rangle, 0,    \langle 1, 2, 5\rangle + \langle 1, 3, 7\rangle + \langle 1, 4, 5\rangle +
\langle 1, 4, 7\rangle + \langle 3, 4, 7\rangle\},&\\
&\{\langle 2, 6\rangle,     0, \langle 1, 2, 6\rangle + \langle 1, 3, 6\rangle + \langle 3, 4, 10\rangle +
 \langle 3, 6, 10\rangle +
      \langle 4, 6, 10\rangle\},&\\
&    \{\langle 2, 8\rangle, 0, \langle 1, 2, 8\rangle\},
\{\langle 2, 11\rangle, 0,    \langle 1, 2, 11\rangle\},&\\
   &   \{\langle 3, 5\rangle, 0,    \langle 1, 3, 5\rangle + \langle 1, 3, 7\rangle+ \langle 1, 4, 5\rangle +
       \langle 1, 4, 7\rangle+ \langle 3, 4, 7\rangle\},&\\
&      \{\langle 3, 6\rangle,     0, \langle 3, 4, 10\rangle + \langle 3, 6, 10\rangle + \langle 4, 6, 10\rangle\},
      \{\langle 3, 7\rangle, 0,    \langle 1, 3, 7\rangle\},&\\
&\{\langle 3, 8\rangle, 0, \langle 1, 3, 8\rangle\},
      \{\langle 3, 11\rangle, 0,    \langle 1, 3, 11\rangle\},&\\
&\{\langle 4, 7\rangle, 0,    \langle 1, 3, 7\rangle + \langle 3, 4, 7\rangle\},
       \{\langle 4, 10\rangle, 0,    \langle 3, 4, 10\rangle\},&\\
&       \{\langle 4, 11\rangle, 0,    \langle 1, 3, 7\rangle + \langle 1, 4, 7\rangle + \langle 1, 4, 11\rangle
+ \langle 3, 4, 7\rangle\},&\\
    &   \{\langle 5, 6\rangle, 0,    \langle 4, 5, 9\rangle + \langle 4, 6, 9\rangle + \langle 5, 6, 9\rangle\},
      \{\langle 5, 9\rangle, 0,    \langle 4, 5, 9\rangle\},&\\
   &    \{\langle 5, 11\rangle, 0,    \langle 1, 3, 7\rangle + \langle 1, 4, 5\rangle
       +\langle 1, 4, 7\rangle + \langle 1, 5, 11\rangle + \langle 3, 4, 7\rangle\},&\\
&       \{\langle 6,9\rangle, 0, \langle 4, 6, 9\rangle\},
 \{\langle 6, 10\rangle, 0, \langle 3, 4, 10\rangle + \langle 4, 6, 10\rangle\},&\\
&       \{\langle 6, 11\rangle, 0,\langle 1, 3, 6\rangle + \langle 1, 6, 11\rangle
       + \langle 3, 4, 10\rangle + \langle 3, 6, 10\rangle +  \langle 4, 6, 10\rangle\},&\\
&      \{\langle 1, 2, 3\rangle, \langle 1, 2, 3\rangle, 0\},
    \{\langle 1, 3, 4\rangle, \langle 1, 3, 4\rangle,    0\},&\\
     & \{\langle 1, 4, 6\rangle, \langle 1, 4, 6\rangle, 0\},
    \{\langle 1, 5, 6\rangle, \langle 1, 5, 6\rangle,    0\},&\\
&    \{\langle 2, 3, 4\rangle,    \langle 1, 2, 3\rangle + \langle 1, 3, 4\rangle,     \langle 1, 2, 3, 4\rangle\},
    \{\langle 2, 3, 5\rangle,    \langle 1, 2, 3\rangle,  \langle 1, 2, 3, 5\rangle\},&\\
&    \{\langle 2, 3, 6\rangle,    \langle 1, 2, 3\rangle,  \langle 1, 2, 3, 6\rangle\},
    \{\langle 2, 3, 11\rangle, \langle 1, 2, 3\rangle,     \langle 1, 2, 3, 11\rangle\},&\\
&  \{\langle 2, 4, 5\rangle, 0,    \langle 1, 2, 4, 5\rangle\},
     \{\langle 2, 4, 6\rangle, \langle 1, 4, 6\rangle,     \langle 1, 2, 4, 6\rangle\},&\\
 &     \{\langle 2, 4, 11\rangle, 0,    \langle 1, 2, 4, 11\rangle\},
   \{\langle 2, 5, 6\rangle, \langle 1, 5, 6\rangle,     \langle 1, 2, 5, 6\rangle\},&\\
     &\{\langle 2, 5, 11\rangle, 0,    \langle 1, 2, 5, 11\rangle\}, 
     \{\langle 2, 6, 11\rangle, 0, \langle 1, 2, 6, 11\rangle\},&\\
   &   \{\langle 3, 4, 5\rangle,    \langle 1, 3, 4\rangle,  \langle 1, 3, 4, 5\rangle\},
     \{\langle 3, 4, 6\rangle,    \langle 1, 3, 4\rangle + \langle 1, 4, 6\rangle,      \langle 1, 3, 4, 6\rangle\},&\\
 &   \{\langle 3, 4, 11\rangle, \langle 1, 3, 4\rangle,    \langle 1, 3, 4, 11\rangle\},
    \{\langle 3, 5, 6\rangle, \langle 1, 5, 6\rangle,    \langle 1, 3, 5, 6\rangle\},&\\
&    \{\langle 3, 5, 11\rangle, 0,   \langle 1, 3, 5, 11\rangle\},
    \{\langle 3, 6, 11\rangle, 0,   \langle 1, 3, 6, 11\rangle\},&\\
    &\{\langle 4, 5, 6\rangle, \langle 1, 4, 6\rangle + \langle 1, 5, 6\rangle, \langle 1, 4, 5, 6\rangle\}, 
     \{\langle 4, 5, 11\rangle, 0, \langle 1, 4, 5, 11\rangle\}, &\\
    &\{\langle 4, 6, 11\rangle,     \langle 1, 4, 6\rangle,    \langle 1, 4, 6, 11\rangle\},
    \{\langle 5, 6, 11\rangle,    \langle 1, 5, 6\rangle, \langle 1, 5, 6, 11\rangle\},&\\
&    \{\langle 2, 3, 4, 5\rangle, 0,    \langle 1, 2, 3, 4, 5\rangle\},
\{\langle 2, 3, 4, 6\rangle, 0,\langle 1, 2, 3, 4, 6\rangle\},&\\
 &   \{\langle 2, 3, 4, 11\rangle,    0, \langle 1, 2, 3, 4, 11\rangle\},
    \{\langle 2, 3, 5, 6\rangle, 0,    \langle 1, 2, 3, 5, 6\rangle\},&\\
   & \{\langle 2, 3, 5, 11\rangle, 0,    \langle 1, 2, 3, 5, 11\rangle\},
    \{\langle 2, 3, 6, 11\rangle, 0,    \langle 1, 2, 3, 6, 11\rangle\},&\\
&   \{\langle 2, 4, 5, 6\rangle, 0, \langle 1, 2, 4, 5, 6\rangle\},
    \{\langle 2, 4, 5, 11\rangle,     0, \langle 1, 2, 4, 5, 11\rangle \},&\\
   & \{\langle 2, 4, 6, 11\rangle, 0,    \langle 1, 2, 4, 6, 11\rangle\},
    \{\langle 2, 5, 6, 11\rangle, 0,    \langle 1, 2, 5, 6, 11\rangle\},&\\
&\{\langle 3, 4, 5, 6\rangle, 0, \langle 1, 3, 4, 5, 6\rangle\},
    \{\langle 3, 4, 5, 11\rangle,     0, \langle 1, 3, 4, 5, 11\rangle\},&\\
&    \{\langle 3, 4, 6, 11\rangle, 0,    \langle 1, 3, 4, 6, 11\rangle\},
    \{\langle 3, 5, 6, 11\rangle, 0,    \langle 1, 3, 5, 6, 11\rangle\},&\\
    &\{\langle 4, 5, 6, 11\rangle, 0,    \langle 1, 4, 5, 6, 11\rangle\},
   \{\langle 2, 3, 4, 5, 6\rangle, 0,    \langle 1, 2, 3, 4, 5, 6\rangle\},&\\
    &\{\langle 2, 3, 4, 5, 11\rangle, 0,    \langle 1, 2, 3, 4, 5, 11\rangle\},
    \{\langle 2, 3, 4, 6, 11\rangle, 0,    \langle 1, 2, 3, 4, 6, 11\rangle\},&\\
  &\{\langle 2, 3, 5, 6, 11\rangle, 0,    \langle 1, 2, 3, 5, 6, 11\rangle\},
    \{\langle 2, 4, 5, 6, 11\rangle, 0,    \langle 1, 2, 4, 5, 6, 11\rangle\},&\\
  & \{\langle 3, 4, 5, 6, 11\rangle, 0,    \langle 1, 3, 4, 5, 6, 11\rangle\},
  \{\langle 2, 3, 4, 5, 6, 11\rangle, \langle 2, 3, 4, 5, 6, 11\rangle,    0\}&\}\,.
\end{eqnarray*}
Notice that if a simplex of $K$ doesn't appear in this list, it is because its image under
$f$ and $\phi$ is null.
The representative cycle of every homology class is:
\begin{eqnarray*}
&&g\langle 1\rangle=\langle 1\rangle\\
&& g\langle 1, 2, 3\rangle= \langle 1, 2, 3\rangle +\langle 1, 2, 8\rangle +
  \langle 1, 3, 8] + \langle 2, 3, 8\rangle\\
&& g\langle 1, 3, 4\rangle=
    \langle 1, 3, 4\rangle + \langle 1, 3, 7\rangle+ \langle 1, 4, 7\rangle + \langle 3, 4, 7\rangle\\
  &&  g\langle 1, 4, 6\rangle=
    \langle 1, 3, 4\rangle + \langle 1, 3, 6\rangle + \langle 1, 4, 6\rangle + \langle 3, 4, 10\rangle +
    \langle 3, 6, 10\rangle +
      \langle 4, 6, 10\rangle\\
&&  g\langle 1, 5, 6\rangle=
    \langle 1, 4, 5\rangle +\langle 1, 4, 6\rangle + \langle 1, 5, 6\rangle+ \langle 4, 5, 9\rangle
    + \langle 4, 6, 9\rangle +
      \langle 5, 6, 9\rangle\\
&& g\langle 2, 3, 4, 5, 6, 11\rangle=
    \langle 1, 2, 3, 4, 5, 6\rangle + \langle 1, 2, 3, 4, 5, 11\rangle + \langle 1, 2, 3, 4, 6, 11\rangle\\
&&\mbox{ } +
      \langle 1, 2, 3, 5, 6, 11\rangle + \langle 1, 2, 4, 5, 6, 11\rangle + \langle 1, 3, 4, 5, 6, 11\rangle
      +
      \langle 2, 3, 4, 5, 6, 11\rangle\,.
     \end{eqnarray*}
  A base of the kernel of $Sq^2H^2K$ is:
$$\{\langle 1, 2, 3\rangle^*, \langle 1, 3, 4\rangle^*, \langle 1, 4, 6\rangle^*,
\langle 1, 5, 6\rangle^*\}\,.$$
Now, given an element $\alpha$ of this kernel, we first have to compute
$c=g^*\alpha$.
Let us study a concrete example with all the details.
Let us take $\alpha=\langle 1, 2, 3\rangle^*+\langle 1, 5, 6\rangle^*$. Then
\begin{eqnarray*}
c=g^*\alpha=\alpha f&=&\langle 1,2,3 \rangle^*+\langle 1,5,6 \rangle^*+\langle 2,3,4 \rangle^*+\langle 2,3,5 \rangle^*
+\langle 2,3,6 \rangle^*\\&&+\langle 2,3,11 \rangle^*+\langle 2,5,6 \rangle^*
+\langle 3,5,6 \rangle^*+\langle 4,5,6 \rangle^*+\langle 5,6,11 \rangle^*\,.\end{eqnarray*}
We now compute the cochains of the 3rd step of the algorithm for computing $\Psi_2$.
$$\delta b=c\smile c=\langle 1,2,3,5,6\rangle^*+\langle 2,3,4,5,6\rangle^*+\langle 2,3,5,6,11\rangle^*$$ $$
 b=(c\smile c)\phi=\langle 2,3,5,6\rangle^*$$
Then, we have that $b\smile_1 b=0$ and $b\smile_2\delta b=0$.
On the other hand, $\delta\eta=c\smile_1 c=0$ therefore $\eta\smile\delta\eta=0$.
We thus get,
$$w=f^*(E_3c^4)=(E_3c^4)g
=\langle 1, 2, 3, 4, 5,6 \rangle^*g=\langle 2,3,4,5,6,11\rangle^*\,.$$
Therefore,
$\Psi_2(\langle 1, 2, 3\rangle^*+\langle 1, 5, 6\rangle^*)=\langle 2,3,4,5,6,11\rangle^*$.
Finally, observe that  since there are no classes of cohomology of dimension $3$, then
$\langle 2,3,4,5,6,11\rangle^*\not\in  \,Im\, Sq^2 H^3K$.

\section{Some Comments}

All these results can be given in a more general framework working not necessarily with finite
simplicial complexes. Nevertheless, a contraction from the
(co)chain complex associated to the simplicial complex to its (co)homology must exist
in order to develop the method.

In this paper, the ground ring is ${\bf Z}_2$ for simplicity, but  the same process can be done working
with any field as the ground ring.
For example, let  ${\bf Z}_p$ ($p$ being a prime) be the group of coefficients.
From the combinatorial formulae for
the reduced $p$th powers $P_i$ \cite{Ste47,SE62}
at cochain level in terms of face operators  established in
\cite{GR99,Gon00} and the algorithm for computing the chain contraction $(f,g,\phi)$ from
$C_*(K;{\bf Z}_p)$ to $H_*(K;{\bf Z}_p)$, Steenrod cohomology operations can effectively
be  computed.   Let $\alpha^*\in H^q(K;{\bf Z}_p)$, for calculating the
cohomology class ${\cal P}_i(\alpha^*)$ with
 $\alpha^*\in H^q(K;{\bf Z}_p)$,
we only have to compute $P_i(\alpha f)g$.

Finally,
in order to obtain the image of
any cohomology operation at cochain level
 over a representative
cocycle using our formulae, we have to compute them on a base of $C_*(K)$
in the desired dimension.
A way of decreasing the complexity of this
is to do a ``topological" thinning of the simplicial complex $K$
in order to obtain a thinned simplicial subcomplex $M_{\mbox{\scsz top}}K$ of $K$ (such that
there exists a chain contraction from $C_*K$ to $C_*(M_{\mbox{\scsz top}}K)$). Two examples of thinning in this way
are edge contractions (example (a)) and simplicial collapses (example (b)).
Therefore, we can apply our machinery to compute cohomology operations in the thinned
simplicial complex $M_{\mbox{\scsz top}}K$ and then, the results can be
easily interpreted in the ``big" simplicial complex $K$ via composition of contractions.

\end{document}